\documentclass[12pt]{article}

\usepackage[cp1251]{inputenc}
\usepackage[english]{babel}
\usepackage{indentfirst}
\usepackage{latexsym}
\usepackage{amssymb}
\usepackage[reqno,fleqn]{amsmath}
\usepackage{amsthm}
\usepackage{enumerate}
\usepackage{authblk}
\usepackage{cite}
\usepackage{appendix}

\newtheorem{theorem}{Theorem}
\newtheorem{corollary}{Corollary}
\newtheorem{lemma}{Lemma}

\newtheorem{proposition}{Proposition}

\begin{document}

\title{The principle of invariance in the Donsker form to the partial sum processes of finite order moving averages}

\author[$,a$]{N. S. Arkashov \thanks{Corresponding author.\\
\hspace*{5mm}E-mail address: nicky1978@mail.ru (N. S. Arkashov).}}
\affil[$a$]{Novosibirsk State Technical University,
Karl Marx Ave., 20, Novosibirsk 630073, Russia.}

\date{}
\maketitle
\leftline{UDC 519.214}
\begin{abstract}
We consider the process of partial sums of moving averages of finite order with a regular varying memory function, constructed from a stationary sequence, variance of the sum of which is a regularly varying function. We study the Gaussian approximation of this process of partial sums with the aid of a certain class of Gaussian processes, and obtain sufficient conditions for the $C$-convergence in the invariance principle in the Donsker form.\medskip
\end{abstract}

\textbf{keywords:} invariance principle, fractional Brownian motion, moving average, Gaussian process, memory function, regular varying function

\section{Introduction and Statement of Main Results}
In the present paper, sufficient conditions are obtained for $C$-convergence in the metric space $D[0,1]$ of a normalized process of partial sums of moving averages of finite order to a Gaussian process constructed from a fractional Brownian motion. We will use constructions, and also adhere to the \cite{Arkashov_Arka} presentation scheme. A significant drawback of the estimates of the proximity of processes obtained in \cite{Arkashov_Arka} in the principle of invariance in the Strassen form is the ``local conditions'' for each element of a non-random sequence that forms two-sided moving averages. Note that, in most cases, it is not possible to restore information about the behavior of this non-random sequence during statistical data analysis. This shortcoming was the motivation for writing this article, in which the mentioned ``local conditions'' are replaced by an ``integral condition'' for the regular behavior of the variance of the process of partial sums of two-sided moving averages. Note that from the assumption that this ``integral condition'' is satisfied, one can obtain estimates for the scale parameter and the Hurst parameter (see \cite{Arkashov_zvm}).

Let $\{X_j;\;j\in\mathbb{Z}\}$ be a stationary (in the narrow sense) sequence of random variables represented as a two-sided moving average:
\begin{equation}
\label{Arkashov_1in} X_j=\sum_{k=-\infty}^{\infty}a_{j-k}\xi_k,
\end{equation}
where $\{\xi_k;\;k\in \mathbb{Z}\}$ is a sequence of i. i. d. random variables with zero means and unit variances, $\{a_k;\;k\in \mathbb{Z}\}$ is a nonrandom square summable sequence of real numbers.

We introduce the notation: $S_n=\sum_{i=1}^n X_i$ for $n\geq 1$, $S_0=0$. In what follows, we assume that ${\bf Var}(S_n)$ is a regularly varying function
\begin{equation}\label{Arkashov_usH}
{\bf Var}(S_n)=h^2(n)n^{2H},
\end{equation}
where $h$ is a slowly varying function at $+\infty$ (e.g., see \cite{Arkashov_Fell}) and $H\in (0,1)$.
\\\\
Denote by $\mathcal{R}_{\nu}$ the class of regularly varying functions with exponent $\nu\geq0$, defined on the interval $[0,+\infty)$, different from constants having the following representation: $l(t)t^{\nu}$, where $l(t)$ is a slowly varying at $+\infty$, non-negative and non-decreasing function on $[0,+\infty)$ (we assume that $0^0=1$).
\\\\

We construct a sequence of moving average of finite order
\begin{equation}\label{Arkashov_sk}
v_0=0,~v_k=\sum_{i=0}^{k-1} X_{k-i}\Delta M(i),~k\geq 1,
\end{equation}
where $M\in \mathcal{R}_{\nu}$ and $\Delta M(i)=M(i+1)-M(i)$. Following \cite{Arkashov_Olem}, \cite{Arkashov_Nigm}, the function $M$ will be called the \emph{memory function}.
We define the partial sums process of the sequence (\ref{Arkashov_sk})
\begin{equation}\label{Arkashov_loma}
R_{n}=\sum_{k=0}^{n}v_k,~n=0,1,\dots.
\end{equation}
Note that a possible ``physical meaning'' for the $\{R_{n}\}$ random walk is presented in \cite{Arkashov_Sel}.

By $B_H(t)$, we denote the fractional Brownian motion, i.e., centered Gaussian process with covariance function
\begin{equation}\label{Arkashov_yop}
R(t, s)=\frac{1}{2}\left(t^{2H}+s^{2H}-|t-s|^{2H}\right),
\end{equation}
where $H\in (0,1)$.

Before stating the main result of the paper, we present the assertion underlying its proof.
\begin{proposition}\label{Arkashov_ccon} 
Let condition {\rm(\ref{Arkashov_usH})} and ${\bf E}|\xi_1|^\alpha<\infty$ be satisfied for some
$\alpha$ such that $\alpha\geq 2$ and $\alpha H>1$.
Then, as $n\to\infty$, the processes $S_{[nt]}/\sqrt{{\bf Var}(S_n)}$ $C$-converge in $D[0,1]$ to $B_H(t)$.
\end{proposition}
%For a real number $x$, let $[x]$ denote the largest integer not exceeding $x$.
In Proposition \ref{Arkashov_ccon} the expression $[x]$ denotes the largest integer not exceeding the real number $x$.

Recall that by $C$-convergence in $D[0,1]$ we mean the weak convergence of distributions of measurable (in Skorokhod's topology) functionals on $D[0,1]$ that are continuous in uniform topology at the points of the space $C[0,1]$ (e.g., see \cite{Arkashov_Sah}).

Note that the specified $C$-convergence is established in \cite{Arkashov_Dav} under more stringent moment constraints than those in Proposition \ref{Arkashov_ccon}.
We also note that in the case $H>1/2$ the assertion of Proposition \ref{Arkashov_ccon} was proved in \cite[Theorem 2]{Arkashov_Kons}.

Let us define the Gaussian process $Z_{\nu,H}$:
$$Z_{\nu,H}(t)= \nu\int_0^t B_H(t-s)s^{\nu-1}\,ds.$$
We introduce the notation
$$r_n(t)=\frac{R_{[nt]}}{h(n)n^H M(n)},~~t\in [0,1],~n=1,2,\dots.$$

\begin{theorem}\label{Arkashov_post}
Let the condition of Proposition {\rm\ref{Arkashov_ccon}} be fulfilled. Then

{\rm 1)} if $M\in \mathcal{R}_{\nu},~\nu>0$, then the processes $r_n(t)$ $C$-converge in $D[0,1]$ to $Z_{\nu,H}(t)$ as $n\to\infty${\rm;}

{\rm 2)} if $M\in \mathcal{R}_{0}$, then the processes $r_n(t)$ $C$-converge in $D[0,1]$ to $\left(1-\frac{M(0)}{M(+\infty)}\right)B_H(t)$ as $n\to\infty$.
\end{theorem}

\begin{corollary}\label{Arkashov_cr1}
Let condition {\rm(\ref{Arkashov_usH})} be fulfilled. Then, as $n\rightarrow\infty$,

{\rm 1)} for the case of $M\in\mathcal{R}_{\nu},~\nu>0$ the following equivalence is valid
$${\bf Var}(R_n)\sim {\bf Var}(Z_{\nu,H}(1)) M^2(n) h^2(n)n^{2H};$$

{\rm 2)} for the case of $M\in \mathcal{R}_{0}$ it holds that
$${\bf Var}(R_n)\sim \left(1-\frac{M(0)}{M(+\infty)}\right)^2 M^2(n)h^2(n)n^{2H}.$$
\end{corollary}

\section{Proofs.}

Let us divide the proof of Proposition \ref{Arkashov_ccon} into lemmas \ref{Arkashov_l10}---\ref{Arkashov_l16}. Note that the proof of this proposition follows the scheme of the proof of Theorem 2 in \cite{Arkashov_Arb}. Introduce the following notations:
$$s_{n}(t)=S_{[nt]}/\sqrt{{\bf Var}(S_n)},~~~ t\in
[0,\,1],~~~n=1,\,2,\,\dots,$$
$$A_{k,n}(t)=(n^H h(n))^{-1}\sum_{j=-k+1}^{-k+[nt]}a_j.$$
Using the above notation, we immediately obtain (see \cite{Arkashov_Kons})
$$s_{n}(t)=\sum_{k\in\mathbb{Z}}A_{k,n}(t)\xi_k.$$

\begin{lemma}
\label{Arkashov_l11} For all $1\geq t\geq\tau\geq0$ it holds that
$$\sum_{k\in\mathbb{Z}}(A_{k,n}(t)-A_{k,n}(\tau))^2=
\frac{([nt]-[n\tau])^{2H}h^2([nt]-[n\tau])}{n^{2H}h^2(n)}.$$
\end{lemma}
\begin{proof}
We have
$$\sum_{k\in\mathbb{Z}}(A_{k,n}(t)-A_{k,n}(\tau))^2=n^{-2H}h^{-2}(n)
\sum_{i\in\mathbb{Z}}\left(\sum_{j=-i+1}^{-i+[nt]-[n\tau]}a_j\right)^2=
\frac{{\bf Var}(S_{[nt]-[n\tau]})}{n^{2H}h^{2}(n)}.$$ Using (\ref{Arkashov_usH}), we obtain the assertion of the lemma.

\end{proof}

\begin{lemma}
\label{Arkashov_l10} For all $1\geq t\geq\tau\geq0$ the following convergence is valid as
$n\rightarrow\infty$
$${\bf E}s_{n}(t)s_{n}(\tau)\rightarrow {\bf E}B_H(t)B_H(\tau).$$
\end{lemma}
\begin{proof} We have
\begin{gather*}
{\bf E}(s_{n}(t)-s_{n}(\tau))^2=\sum_{k\in \mathbb{Z}}(A_{k,n}(t)-A_{k,n}(\tau))^2\\
=\frac{([nt]-[n\tau])^{2H}h^2([nt]-[n\tau])}{n^{2H}h^2(n)},
\end{gather*}
where the last equality follows from Lemma \ref{Arkashov_l11}.
From here we get
$${\bf E}(s_{n}(t)-s_{n}(\tau))^2\rightarrow(t-\tau)^{2H}=
{\bf E}(B_H(t)-B_H(\tau))^2.$$
Using the convergence of the second moments of the one-dimensional projections of the random processes $s_{n}(t)$ we deduce the equality
\begin{equation*}
\begin{split}
2{\bf E}s_{n}(t)s_{n}(\tau)&={\bf E}(s_{n}(t))^2+{\bf
E}(s_{n}(\tau))^2\\
&-{\bf E}(s_{n}(t)-s_{n}(\tau))^2\rightarrow 2{\bf
E}B_H(t)B_H(\tau).
\end{split}
\end{equation*}

\end{proof}

\begin{lemma}[{\cite{Arkashov_Petr}}]\label{Arkashov_l12}
Let $\{X_k\}_{k=1,\dots,n}$ be independent centered random variables and
$E|X_k|^\alpha <+\infty$ for all $k$ and some $\alpha\geq2$. Put
$$S_n=\sum_{k=1}^nX_k,\;
M_{\alpha,n}=\sum_{k=1}^n {\bf
E}|X_k|^\alpha,\;B_n=\sum_{k=1}^n{\bf E}X_k^2.$$ Then
$${\bf E}|S_n|^\alpha\leq c(\alpha)(M_{\alpha,n}+B_n^{\alpha/2}),$$ where
$c(\alpha)$ is a positive constant depending only on
$\alpha$.
\end{lemma}
\begin{lemma}
\label{Arkashov_l13}  The following inequality is valid:
$${\bf E}|s_{n}(t)-s_{n}(\tau)|^\alpha\leq
C\left(\frac{[nt]-[n\tau]}{n}\right)^{\alpha H}
\left(\frac{h([nt]-[n\tau])}{h(n)}\right)^{\alpha},$$ where $C$
is a constant depending on the distributions of $\xi_0$, $\alpha$ and
$\{a_i\}$.
\end{lemma}
\begin{proof}
By Lemma \ref{Arkashov_l12} and the Fatou theorem we obtain
\begin{equation}\label{Arkashov_gavv}
\begin{split}
&{\bf E}|\sum_{k\in\mathbb{Z}}(A_{k,n}(t)-A_{k,n}(\tau))\xi_k|^\alpha\\
&\leq c(\alpha)\left(\sum_{k\in\mathbb{Z}}|A_{k,n}(t)-A_{k,n}(\tau)|^\alpha
{\bf E}|\xi_0|^\alpha+
\left(\sum_{k\in\mathbb{Z}}(A_{k,n}(t)-A_{k,n}(\tau))^2\right)^{\alpha/2}\right)\\
&\leq c(\alpha)(1+{\bf E}|\xi_0|^\alpha)
\left(\sum_{k\in\mathbb{Z}}(A_{k,n}(t)-A_{k,n}(\tau))^2\right)^{\alpha/2}.
\end{split}
\end{equation}
Note that we applied here the following elementary inequality:
$$\sum_{k\in\mathbb{Z}}|b_k|^\gamma\leq \left(\sum_{k\in\mathbb{Z}}|b_k|\right)^\gamma,\;\gamma\geq1.$$
Next, applying Lemma \ref{Arkashov_l11} to the right-hand side of the last inequality in (\ref{Arkashov_gavv}), we obtain the assertion of the lemma.

\end{proof}

\begin{lemma}\label{Arkashov_mmfu}
Let $g$ be a slowly varying function defined on $[1,+\infty)$ and $\varepsilon$ be an arbitrary positive number. Then
the sequence $$\left\{\frac{\max_{1\leq k \leq
n}g(k)k^{\varepsilon}}{g(n)n^{\varepsilon}}\right\}$$ is bounded
for all $n\geq 1$.
\end{lemma}
\begin{proof}
Using Statement 4 from {\cite[section 1.5]{Arkashov_Seneta}}, we find that there exists $B$ such that $\frac{\sup_{B\leq x \leq
n}g(x)x^{\varepsilon}}{g(n)n^{\varepsilon}}\rightarrow 1$ holds as $n\rightarrow +\infty$, which immediately implies the assertion of the lemma.
\end{proof}

From Lemma \ref{Arkashov_l13} and \ref{Arkashov_mmfu} we immediately obtain the following assertion.

\begin{lemma}\label{Arkashov_l14}
Let $\alpha H>1$. Then the following inequality is valid
$${\bf E}|s_{n}(t)-s_{n}(\tau)|^\alpha\leq
C \left(\frac{[nt]-[n\tau]}{n}\right)^{\frac{\alpha H+1}{2}},$$
where $C$ is a constant depending on the distribution of $\xi_0$, $\alpha$ and $\{a_i\}$.
\end{lemma}
From Lemma \ref{Arkashov_l14}, applying Theorem 3 from \cite[p. 52]{Arkashov_Sah} (where it is necessary to take $[nt]/n$ as the step function $F_n(t)$), we deduce that for any $\varepsilon>0$ it holds
\begin{equation}\label{Arkashov_pltn}
\lim_{\delta\rightarrow 0}\limsup_{n\rightarrow+\infty}{\bf P}\left(\sup_{|t-s|<\delta,~t,s\in
[0,1]}|s_n(t)-s_n(s)|>\varepsilon\right)=0.
\end{equation}

\begin{lemma}[{\cite{Arkashov_Ibr}}]\label{Arkashov_l9}
For all $k\in\mathbb{Z}$ the following is valid
$$|a_{k+1}+...+a_{k+n}|\leq
\left(4\sqrt{{\bf Var}(S_n)}\sum_{k\in\mathbb{Z}}|a_k|^2\left(1+\frac{1}{2\sqrt{{\bf Var}(S_n)}}\right)\right)^{1/2}.$$
\end{lemma}
\begin{lemma}[{\cite[Lemma 7]{Arkashov_Kons}}]
\label{Arkashov_l15} Let $\{b_{ni};\;n\geq1,\;i\in\mathbb{Z}\}$ be an array of real numbers, and let $\{\zeta_{ni};\;n\geq1,\;i\in\mathbb{Z}\}$ be an
array of random variables satisfying the following conditions{\rm:}
\begin{enumerate}
\item[\rm L1.]
$\lim_{n\rightarrow\infty}\sum_{i\in\mathbb{Z}}b_{ni}^2=1${\rm ;}

\item[\rm L2.]
$\lim_{n\rightarrow\infty}\sup_{i\in\mathbb{Z}}|b_{ni}|=0${\rm ;}

\item[\rm L3.] For every $n\geq1$ the sequence
$\{\zeta_{ni},\;i\in\mathbb{Z}\}$ consists of i.i.d. random variables with mean zero and
variance $1${\rm ;}
%${\bf E}\zeta_{n0}=0$ и ${\bf E}\zeta_{n0}^2=1$

\item[\rm L4.] $\lim_{K\rightarrow\infty}\sup_{n\geq1}{\bf
E}\zeta_{n0}^2 {\rm I}(|\zeta_{n0}|>K)=0$.
\end{enumerate}

Then the sums $\sum_{i\in\mathbb{Z}}b_{ni}\zeta_{ni}$
converge in distribution to a standard Gaussian random variable as $n\rightarrow +\infty$.
\end{lemma}

\begin{lemma}\label{Arkashov_l16}
The finite-dimensional distributions of random processes $s_{n}(t)$ converge to the corresponding finite-dimensional distributions of random process $B_H(t)$ as $n\rightarrow\infty$.
\end{lemma}
\begin{proof}
It suffices to prove that $\sum_{i=1}^l c_i s_{n}(t_i)$ converges in distribution to $\sum_{i=1}^l c_iB_H(t_i)$ for every
finite set of numbers $\{c_i;\,i=1,\dots,l\}$. So, we observe first that
$z_n\equiv\sum_{i=1}^lc_is_{n}(t_i)=\sum_{k\in\mathbb{Z}}\sum_{i=1}^l
c_iA_{k,n}(t_i)\xi_k.$ Further (see Lemma \ref{Arkashov_l10})

\begin{equation*}
\begin{split}
{\bf E}z_n^2=\sum_{i,j=1}^lc_ic_j{\bf
E}(s_{n}(t_i)s_{n}(t_j))\rightarrow
&\sum_{i,j=1}^lc_ic_j{\bf E}(B_H(t_i)B_H(t_j))\\
&={\bf E}\left(\sum_{i=1}^l c_i
B_H(t_i)\right)^2,\;n\rightarrow\infty.
\end{split}
\end{equation*}
Use the notation
$$\delta^2={\bf
E}\left(\sum_{i=1}^l c_i
B_H(t_i)\right)^2,\;z_{1n}=z_n/\delta.$$ Thus
%\begin{equation*}
$$\label{Arkashov_1th2} Ez_{1n}^2=\sum_{k\in\mathbb{Z}}\left(\sum_{i=1}^lc_i/\delta
A_{k,n}(t_i)\right)^2\rightarrow1 ,\;n\rightarrow\infty.$$
%\end{equation*}
Now, under the conditions of Lemma \ref{Arkashov_l15} put
$b_{nk}=\sum_{i=1}^lc_i/\delta A_{k,n}(t_i),\;\zeta_{nk}=\xi_k$.
It is easy to see that, in this case, condition $L1$ is valid.
Using Lemma \ref{Arkashov_l9} we find, as $n\rightarrow +\infty$, that
$$\sup_{k}|\sum_{i=1}^l(c_iA_{k,n}(t_i))|=O(n^{-H/2}h^{-1/2}(n)),\;n\rightarrow\infty.$$
Hence, condition ${L2}$ is valid.
It is easy to verify that conditions ${L3}$ and ${ L4}$ hold as well.
So, the random variables $z_n$ converge in distribution to a normal random variable with mean zero
and variance $\delta^2$.

\end{proof}

\begin{proof}[Proof of Proposition {\rm\ref{Arkashov_ccon}}.]
Lemmas \ref{Arkashov_l14}, \ref{Arkashov_l16} and relation (\ref{Arkashov_pltn}) immediately imply the assertion of the proposition.
\end{proof}

Let's proceed to the proof of Theorem \ref{Arkashov_post}. By changing the order of summation in (\ref{Arkashov_loma}), represent $R_n$ as $\sum_{i=0}^n S_{n-i}\Delta M(i)$. From where we derive the representation of $R_n$ as an integral $\int_0^n S_{n-[s]}\,dM(s)$. Further, since $r_n(t)=\frac{R_{[nt]}}{h(n)n^H M(n)}$ we get the representation of $r_n(t)$ as:
\begin{equation}\label{Arkashov_loma1}
\int_0^{[nt]/n} s_n(([nt]-[nu])/n)\,\frac{dM(un)}{M(n)}.
\end{equation}

\begin{lemma}\label{Arkashov_plot}
Let $M\in \mathcal{R}_{\nu},~\nu\geq 0$.
Then, the family of distributions of the stochastic processes $\{r_n(t)\}$ is dense in $D[0,1]$ with respect to the uniform metric.
\end{lemma}
\begin{proof}
With the representation (\ref{Arkashov_loma1}) in mind, we find an upper bound for
\begin{equation}\label{Arkashov_voc2}
\begin{split}
\sup_{|t-s|<\delta,~t,s\in
[0,1]}|\int_0^{[nt]/n}s_n&(([nt]-[nu])/n)\,\frac{dM(un)}{M(n)}\\
&-\int_0^{[ns]/n}s_n(([ns]-[nu])/n)\,\frac{dM(un)}{M(n)}|.
\end{split}
\end{equation}

Let, for definiteness, $s<t$. The obvious inequalities hold
\begin{gather*}
%\begin{split}
\left|\int_0^{[nt]/n}s_n(([nt]-[nu])/n)\,\frac{dM(un)}{M(n)}-\int_0^{[ns]/n}s_n(([ns]-[nu])/n)\,\frac{dM(un)}{M(n)}\right|\\
\leq \int_0^{[ns]/n}|s_n(([nt]-[nu])/n)-s_n(([ns]-[nu])/n)|\,\frac{dM(un)}{M(n)}\\
+\int_{[ns]/n}^{[nt]/n}|s_n(([nt]-[nu])/n)|\,\frac{dM(un)}{M(n)}.
%\end{split}
\end{gather*}

If $|t-s|<\delta$, then, for all sufficiently large $n$, it holds that
$$
\frac{[nt]-[nu]}{n}-\frac{[ns]-[nu]}{n}<2\delta.
$$
Since $u\geq [ns]/n$, we have
$$
\frac{[nt]-[nu]}{n}\leq \frac{[nt]-[ns]}{n}< 2\delta.
$$
Finally, we obtain the following upper bound on (\ref{Arkashov_voc2}):
\begin{gather*}
%\begin{split}
\sup_{|t-s|<\delta,~t,s\in
[0,1]}\left|\int_0^{[nt]/n}s_n(([nt]-[nu])/n)\,\frac{dM(un)}{M(n)}-\int_0^{[ns]/n}s_n(([ns]-[nu])/n)\,\frac{dM(un)}{M(n)}\right|\\
\leq 2\sup_{|t-s|<2\delta,~t,s\in
[0,1]}|s_n(t)-s_n(s)|.
%\end{split}
\end{gather*}
Whence, using the relation (\ref{Arkashov_pltn}), we obtain the assertion of the lemma.
\end{proof}

To prove Theorem \ref{Arkashov_post} (taking into account Lemma \ref{Arkashov_plot}), it is sufficient to prove the weak convergence of the finite-dimensional distributions of the processes $r_n(t)$.

We decompose the proof of convergence of the finite-dimensional distributions in the case $\nu>0$ into Lemmas \ref{Arkashov_sh}--\ref{Arkashov_cons1} and, respectively, in the case $\nu=0$ into Lemmas \ref{Arkashov_sqvle}--\ref{Arkashov_cons2}.

In the proofs of the following propositions, we will use the following simple fact. Since $l$ is a slowly varying at $+\infty$ function, there exists an infinitely small sequence $\{\varepsilon_n\}$ such that
$$
n\varepsilon_n\rightarrow +\infty,~~ l(n\varepsilon_n)/l(n)\rightarrow 1.
$$

Recall that $M\in \mathcal{R}_{\nu}$ for $~\nu>0$ can be represented as $M(t)=t^{\nu}l(t)$, where $l$ is a function that is nondecreasing on the interval $[0,+\infty)$ and slowly varying at $+\infty$.
\begin{lemma}\label{Arkashov_sh}
As $n\rightarrow +\infty$, for every $t\in[0,1]$ it holds that{\rm:}
$$
r_n(t)-\int_0^{[nt]/n}s_n(([nt]-[nu])/n)\,du^{\nu}\stackrel{P}{\rightarrow} 0.
$$
\end{lemma}
\begin{proof}

Let $t\in(0,1]$. The following inequalities are valid
\begin{equation}\label{Arkashov_nrvv}
\begin{split}
&\left|\int_0^{[nt]/n}s_n(([nt]-[nu])/n)u^{\nu}\,\frac{dl(un)}{l(n)}\right|\\
&\leq \int_{\varepsilon_n}^{[nt]/n}|s_n(([nt]-[nu])/n)|u^{\nu}\,\frac{dl(un)}{l(n)}\\
&+\int_0^{\varepsilon_n}|s_n(([nt]-[nu])/n)|u^{\nu}\,\frac{dl(un)}{l(n)}\\
&\leq \sup_{t\in[0,1]}|s_n(t)|\left(\frac{l([nt])}{l(n)}-\frac{l(n\varepsilon_n)}{l(n)}\right)+\sup_{t\in[0,1]}|s_n(t)|\varepsilon_n^{\nu}\stackrel{P}{\rightarrow} 0,
\end{split}
\end{equation}
where the convergence to $0$ is implied by the fact that
$$\sup_{t\in[0,1]}|s_n(t)|\stackrel{d}{\rightarrow}\sup_{t\in[0,1]}|B_H(t)|$$ and
$$\frac{l([nt])}{l(n)}-\frac{l(n\varepsilon_n)}{l(n)}\rightarrow 0.$$

We have
\begin{equation}\label{Arkashov_nrvv1}
\begin{split}
&\int_0^{[nt]/n}s_n(([nt]-[nu])/n)\frac{l(un)}{l(n)}\,du^{\nu}\\
&=\int_0^{[nt]/n}s_n(([nt]-[nu])/n)\left(\frac{l(un)}{l(n)}-1\right)\,du^{\nu}\\
&+\int_0^{[nt]/n}s_n(([nt]-[nu])/n)\,du^{\nu}.
\end{split}
\end{equation}
Let us estimate the first term on the right side (\ref{Arkashov_nrvv1}):
\begin{equation}\label{Arkashov_nrvv2}
\begin{split}
&\left|\int_0^{[nt]/n}s_n(([nt]-[nu])/n)\left(\frac{l(un)}{l(n)}-1\right)\,du^{\nu}\right|\\
&\leq \int_0^{\varepsilon_n}|s_n(([nt]-[nu])/n)|\left(1-\frac{l(un)}{l(n)}\right)\,du^{\nu}\\
&+\int_{\varepsilon_n}^{[nt]/n}|s_n(([nt]-[nu])/n)|\left(1-\frac{l(un)}{l(n)}\right)\,du^{\nu}\\
&\leq \sup_{t\in[0,1]}|s_n(t)|\varepsilon_n^{\nu}+\sup_{t\in[0,1]}|s_n(t)|\left(1-\frac{l(n\varepsilon_n)}{l(n)}\right)\stackrel{P}{\rightarrow} 0.
\end{split}
\end{equation}
Next, we note that
\begin{gather*}
\left|\int_0^{[nt]/n}s_n(([nt]-[nu])/n)\,\frac{dM(un)}{M(n)}-\int_0^{[nt]/n}s_n(([nt]-[nu])/n)\,du^{\nu}\right|\\
\leq\left|\int_0^{[nt]/n}s_n(([nt]-[nu])/n)u^{\nu}\,\frac{dl(un)}{l(n)}\right|+\left|\int_0^{[nt]/n}s_n(([nt]-[nu])/n)\left(\frac{l(un)}{l(n)}-1\right)\,du^{\nu}\right|.
\end{gather*}
Whence, using the relations (\ref{Arkashov_nrvv}) and (\ref{Arkashov_nrvv2}), we obtain the assertion of the lemma.
\end{proof}
Introduce the notation:
$$s^{(1)}_{n}(t)=S_{[nt]+1}/\sqrt{{\bf Var}(S_n)},~~~ t\in
[0,\,1],~~~n=1,\,2,\,\dots.$$

\begin{lemma}\label{Arkashov_eqvi}
Let ${\bf E}|\xi_1|^{\alpha}<+\infty$ for some $\alpha$ such that $\alpha\geq 2$ and $\alpha H>1$. Then, as $n\rightarrow +\infty$, it holds that{\rm:}
$$
\sup_{t\in[0,1]}|s^{(1)}_n(t)-s_n(t)|\stackrel{P}{\rightarrow} 0.
$$
\end{lemma}
\begin{proof}
Consider arbitrary $\varepsilon>0$. The following obvious equality holds:
\begin{equation}\label{Arkashov_ochn}
{\bf P}(\sup_{t\in[0,1]}|s^{(1)}_n(t)-s_n(t)|>\varepsilon)\leq n {\bf P}(|X_1|>\varepsilon h(n)n^H).
\end{equation}
Applying Chebyshev's inequality to the right side (\ref{Arkashov_ochn}), we derive
$$
{\bf P}(\sup_{t\in[0,1]}|s^{(1)}_n(t)-s_n(t)|>\varepsilon)\leq \frac{n^{1-\alpha H}}{\varepsilon^{\alpha} h^{\alpha}(n)}{\bf E}|X_1|^{\alpha}\rightarrow 0.
$$
Note that the finiteness of ${\bf E}|X_1|^{\alpha}$ follows from Lemma \ref{Arkashov_l13}.

\end{proof}

\begin{lemma}\label{Arkashov_sh1}
As $n\rightarrow +\infty$, for every $t\in [0,1]$ it holds that{\rm:}
$$
\left|\int_0^{[nt]/n}s_n(([nt]-[nu])/n)\,du^{\nu}-\int_0^t (-s^{(1)}_{n}(u))\,d(t-u)^{\nu}\right|\stackrel{P}{\rightarrow} 0.
$$
\end{lemma}
\begin{proof}
Let $t\in (0,1]$. Let us represent $\int_0^t s^{(1)}_n(u)\,d(t-u)^{\nu}$ as a sum:
\begin{equation}\label{Arkashov_vbnm}
\begin{split}
\int_0^t s^{(1)}_n(u)\,d(t-u)^{\nu}=&\int_0^{[nt]/n} s^{(1)}_n(u)\,d(t-u)^{\nu}\\
&+\int_{[nt]/n}^t s^{(1)}_n(u)\,d(t-u)^{\nu}.
\end{split}
\end{equation}

We estimate the second term on the right-hand side of (\ref{Arkashov_vbnm}).
Using Lemma \ref{Arkashov_eqvi} and Proposition \ref{Arkashov_ccon}, we have
\begin{equation}\label{Arkashov_vbnm1}
\left|\int_{[nt]/n}^t s^{(1)}_n(u)\,d(t-u)^{\nu}\right|\leq |s^{(1)}_n(t)|(t-[nt]/n)^{\nu}\stackrel{P}{\rightarrow} 0.
\end{equation}
Further, we note that the following equality is valid
\begin{equation}\label{Arkashov_tozh}
\int_0^{[nt]/n}s_n(([nt]-[nu])/n)\,du^{\nu}=-\int_0^{[nt]/n} s^{(1)}_n(u)\,d([nt]/n-u)^{\nu}.
\end{equation}
Using (\ref{Arkashov_tozh}), we conclude that
\begin{equation}\label{Arkashov_vbnm2}
\begin{split}
&\left|\int_0^{[nt]/n}s_n(([nt]-[nu])/n)\,du^{\nu}-\int_0^{[nt]/n} (-s^{(1)}_n(u))\,d(t-u)^{\nu}\right|\\
&= |\int_0^{[nt]/n}s^{(1)}_n(u)\,d((t-u)^{\nu}-([nt]/n-u)^{\nu})|\\
&\leq \sup_{t\in[0,1]}|s^{(1)}_n(t)||(t-[nt]/n)^{\nu}-(t^{\nu}-([nt]/n)^{\nu})|\stackrel{P}{\rightarrow} 0.
\end{split}
\end{equation}
To estimate the third integral in (\ref{Arkashov_vbnm2}) we use the fact that the function $f(u)=(t-u)^{\nu}-([nt]/n-u)^{\nu}$ is monotone on the interval $[0,[nt]/n]$.

It's obvious that
\begin{equation}\label{Arkashov_vbnm3}
\begin{split}
&\left|\int_0^{[nt]/n}s_n(([nt]-[nu])/n)\,du^{\nu}-\int_0^t (-s^{(1)}_{n}(u))\,d(t-u)^{\nu}\right|\\
&\leq |\int_{[nt]/n}^t s^{(1)}_n(u)\,d(t-u)^{\nu}|\\
&+|\int_0^{[nt]/n}s_n(([nt]-[nu])/n)\,du^{\nu}-\int_0^{[nt]/n} (-s^{(1)}_n(u))\,d(t-u)^{\nu}|.
\end{split}
\end{equation}
Applying (\ref{Arkashov_vbnm1}) and (\ref{Arkashov_vbnm2}) to the right-hand side of (\ref{Arkashov_vbnm3}), we obtain the assertion of the lemma.
\end{proof}

\begin{lemma}\label{Arkashov_cons}
Let $\{c_i;\,i=1,\dots,l\}$ be an arbitrary set of real numbers and $\{t_i;\,i=1,\dots,l\}$ be an arbitrary set in the interval $[0,1]$. Then, as $n\rightarrow +\infty$, it holds that
$$\sum_{i=1}^l c_i \int_0^{t_i} s^{(1)}_n(u)\,d(t_i-u)^{\nu}\stackrel{d}{\rightarrow}\sum_{i=1}^l c_i\int_0^{t_i} B_H(u)\,d(t_i-u)^{\nu}.$$
\end{lemma}
\begin{proof}

Define the functional $F:D[0,1]\rightarrow \mathbb{R}$ as follows:
$$
F(f)=\sum_{i=1}^l c_i \int_0^{t_i} f(u)\,d(t_i-u)^{\nu}.
$$
It is clear that $F$ is continuous in the uniform topology at the points of $C[0,1]$ and measurable on $D[0,1]$ in Skorokhod's topology. Therefore, taking into account Lemma \ref{Arkashov_eqvi} and the conditions of Proposition \ref{Arkashov_ccon}, we conclude that $F(s^{(1)}_n)$ converges in distribution to $F(B_H)$ (e.g., see \cite{Arkashov_Sah}), and this immediately implies the assertion of the lemma.
\end{proof}

\begin{lemma}\label{Arkashov_cons1}
As $n\rightarrow\infty$, the finite-dimensional distributions of the stochastic processes $r_n(t)$ converge to the finite-dimensional distributions of process $Z_{\nu,H}(t)$.
\end{lemma}

\begin{proof}
Let $\{c_i;\,i=1,\dots,l\}$ be an arbitrary set of real numbers and $\{t_i;\,i=1,\dots,l\}\subseteq [0,1]$. Lemmas (\ref{Arkashov_sh}) and (\ref{Arkashov_sh1}) imply that
$$
\left|\sum_{i=1}^l c_i r_n(t_i)-\sum_{i=1}^l c_i \int_0^{t_i} (-s^{(1)}_n(u))\,d(t_i-u)^{\nu}\right|\stackrel{P}{\rightarrow}0.
$$
Using this fact, Lemma \ref{Arkashov_cons}, and the well-known Cramer--Wold device (e.g., see \cite{Arkashov_Bill}), we obtain the assertion of the lemma.
\end{proof}

%---------------------------------------------------------------------------------------------------------------------------------------------------------------

Now we proceed to the proof of Theorem \ref{Arkashov_post} for the case $\nu=0$. Recall that in this case the function $M$ is non-negative, non-decreasing, non-constant on the interval $[0;+\infty)$ and slowly varying at $+\infty$.

\begin{lemma}\label{Arkashov_sqvle}
There is a positive constant $C$ such that for any $k,l\in [0,n]$ the following inequality is valid{\rm:}
$$
{\bf E}|s_n(k/n)-s_n(l/n)|^2\leq C\frac{|k-l|^H}{n^H}.
$$
\end{lemma}
\begin{proof}
First of all, we have the equality (see Lemma \ref{Arkashov_l11})
$$
{\bf E}|s_n(k/n)-s_n(l/n)|^2=\frac{|k-l|^{2H} h^2(|k-l|)}{n^{2H}h^2(n)}.
$$

We represent this equality in the form
\begin{equation}\label{Arkashov_sqv}
{\bf E}|s_n(k/n)-s_n(l/n)|^2=\frac{|k-l|^{H} }{n^{H}} \frac{|k-l|^{H} h^2(|k-l|)}{n^{H}h^2(n)}.
\end{equation}
Applying Lemma \ref{Arkashov_mmfu} to the second factor on the right-hand side (\ref{Arkashov_sqv}), we obtain the assertion of the lemma.

\end{proof}

\begin{lemma}\label{Arkashov_diff}
For every $t\in[0,1]$, as $n\rightarrow +\infty$, it holds that{\rm:}
$$
r_n(t)-\left(1-\frac{M(0)}{M(+\infty)}\right)s_n(t)\stackrel{P}{\rightarrow}0.
$$
\end{lemma}
\begin{proof}
Let $t\in(0,1]$. Using the representation (\ref{Arkashov_loma1}) for $r_n(t)$, we have
\begin{equation}\label{Arkashov_rave}
\begin{split}
r_n(t)=&\int_{\varepsilon_n}^{[nt]/n}s_n(([nt]-[nu])/n)\,\frac{dM(un)}{M(n)}\\
&+\int_0^{\varepsilon_n}s_n(([nt]-[nu])/n)\,\frac{dM(un)}{M(n)}.
\end{split}
\end{equation}
Let us estimate the first term in (\ref{Arkashov_rave}):
\begin{equation}\label{Arkashov_rave1}
\begin{split}
&\left|\int_{\varepsilon_n}^{[nt]/n}s_n(([nt]-[nu])/n)\,\frac{dM(un)}{M(n)}\right|\\
&\leq \sup_{t\in[0,1]}|s_n(t)|\left(\frac{M([nt])}{M(n)}-\frac{M(\varepsilon_n n)}{M(n)}\right)
\stackrel{P}{\rightarrow}0.
\end{split}
\end{equation}
We represent the second term in (\ref{Arkashov_rave}) as the sum:
\begin{equation}\label{Arkashov_rave2}
\begin{split}
&\int_0^{\varepsilon_n}s_n(([nt]-[nu])/n)\,\frac{dM(un)}{M(n)}\\
&=\int_0^{\varepsilon_n}(s_n(([nt]-[nu])/n)-s_n(t))\,\frac{dM(un)}{M(n)}+\int_0^{\varepsilon_n}s_n(t)\,\frac{dM(un)}{M(n)}.
\end{split}
\end{equation}
Consider the first term on the right-hand side of (\ref{Arkashov_rave2}).
\begin{equation}\label{Arkashov_rave3}
\begin{split}
&{\bf E}\left(\int_0^{\varepsilon_n}(s_n(([nt]-[nu])/n)-s_n(t))\,\frac{dM(un)}{M(n)}\right)^2\\
&=\int_0^{\varepsilon_n}\int_0^{\varepsilon_n}{\bf E}s_n(t,u)s_n(t,v) \,\frac{dM(un)}{M(n)}\frac{dM(vn)}{M(n)},
\end{split}
\end{equation}
where $s_n(t,u)=s_n(([nt]-[nu])/n)-s_n(t)$. Further, using the obvious inequality: ${\bf E}s_n(t,u)s_n(t,v)\leq \frac{1}{2}({\bf E}s^2_n(t,u)+{ \bf E}s^2_n(t,v))$, as well as Lemma \ref{Arkashov_sqvle}, we obtain ${\bf E}s_n(t,u)s_n(t,v)\leq \frac{ C}{2}(u^H+v^H)$. Applying the last inequality to the right side (\ref{Arkashov_rave3}), we get
\begin{equation}\label{Arkashov_rave4}
\begin{split}
&{\bf E}\left(\int_0^{\varepsilon_n}(s_n(([nt]-[nu])/n)-s_n(t))\,\frac{dM(un)}{M(n)}\right)^2\\
&\leq\frac{C}{2}\int_0^{\varepsilon_n}\int_0^{\varepsilon_n} (u^H+v^H)\,\frac{dM(un)}{M(n)}\frac{dM(vn)}{M(n)}\\
&\leq C\varepsilon_n^H\left(\frac{M(\varepsilon_n n)}{M(n)}-\frac{M(0)}{M(n)}\right)^2 \rightarrow 0,~n\rightarrow +\infty.
\end{split}
\end{equation}
From (\ref{Arkashov_rave4}) it follows that as $n\rightarrow +\infty$
\begin{equation}\label{Arkashov_rave41}
\int_0^{\varepsilon_n}(s_n(([nt]-[nu])/n)-s_n(t))\,\frac{dM(un)}{M(n)}\stackrel{P}{\rightarrow}0.
\end{equation}
Consider the third term on the right side (\ref{Arkashov_rave2}). It's obvious that
\begin{equation}\label{Arkashov_rave5}
\begin{split}
&\int_0^{\varepsilon_n}s_n(t)\,\frac{dM(un)}{M(n)}-s_n(t)\left(1-\frac{M(0)}{M(+\infty)}\right)\\
&=s_n(t)\left(\frac{M(\varepsilon_n n)}{M(n)}-\frac{M(0)}{M(n)}\right)-s_n(t)\left(1-\frac{M(0)}{M(+\infty)}\right)\stackrel{P}{\rightarrow}0.
\end{split}
\end{equation}

Applying the relations (\ref{Arkashov_rave41}), (\ref{Arkashov_rave5}) to (\ref{Arkashov_rave2}) we derive
\begin{equation}\label{Arkashov_nwv}
\int_0^{\varepsilon_n}s_n(([nt]-[nu])/n)\,\frac{dM(un)}{M(n)}-s_n(t)\left(1-\frac{M(0)}{M(+\infty)}\right)\stackrel{P}{\rightarrow}0.
\end{equation}

Next, applying the relations (\ref{Arkashov_rave1}), (\ref{Arkashov_nwv}) to (\ref{Arkashov_rave}), we obtain the assertion of the lemma.

\end{proof}

\begin{lemma} \label{Arkashov_cons2}
The finite-dimensional distributions of the stochastic processes $r_n(t)$ converge to the finite-dimensional distributions of the process $\left(1-\frac{M(0)}{M(+\infty)}\right)B_H(t)$ as $n\rightarrow +\infty$.
\end{lemma}
\begin{proof}
Let $\{c_i;\,i=1,\dots, l\}$ be an arbitrary set of real numbers and $\{t_i;\,i=1,\dots, l\}$ be an arbitrary set in the interval $[0,1]$. Lemma \ref{Arkashov_diff} implies that
$$
\sum_{i=1}^l c_i r_n(t_i)-\sum_{i=1}^l c_i s_n(t_i)\left(1-\frac{M(0)}{M(+\infty)}\right)\stackrel{P}{\rightarrow}0.
$$
Whence and from Lemma \ref{Arkashov_l16}, as well as the Cramer--Wold device, the assertion of the lemma follows.
\end{proof}

\begin{proof}[Proof of Corollary {\rm\ref{Arkashov_cr1}}.]
First, we note that the assertion of Theorem \ref{Arkashov_post} holds for Gaussian analogs of the processes $r_n(t)$
when the random variables $\xi_k$ in (\ref{Arkashov_1in}) are the standard Gaussian; in this case, ${\bf Var}(R_n)$ does not change.

If $M\in \mathcal{R}_{\nu},~\nu>0$,
then, as $n\to\infty$, we have convergence in distribution of the sequence $r_n(1)$ to $Z_{\nu,H}(1)$. Moreover, if $M\in \mathcal{R}_{0}$, then the sequence $r_n(1)$ converge in distribution to $\left(1-\frac{M(0)}{M(+\infty)}\right)B_H(1)$ (recall that $r_n(1)=\frac{R_{n}}{h(n)n^H M(n)}$).

Next, we note that the random variable
$R_n/\sqrt{{\bf Var}(R_n)}\sim {\cal N}(0,1)$ for all $n$. Finally, we obtain the equivalence of the normalizing coefficients $R_n$:
$${\bf Var}(R_n)\sim \sigma^2 h^2(n)n^{2H} M^2(n),\;n\rightarrow\infty,$$
where
\begin{equation*}{\label{L}}
\sigma^2=\left\{
 \begin{array}{rl}
{\bf Var}(Z_{\nu,H}(1)),&~~\mbox{if }~~~~M\in \mathcal{R}_{\nu},~\nu>0,\\
(1-M(0)/M(+\infty))^2,&~~\mbox{if }~~~~M\in \mathcal{R}_{0}.
  \end{array} \right.
\end{equation*}

\end{proof}

\end{document}